\documentclass[11pt,reqno]{amsart}

\theoremstyle{plain}
\newtheorem{theorem} {Theorem}[section]
\newtheorem{lemma}[theorem] {Lemma}

\newtheorem{corollary}[theorem] {Corollary}

\theoremstyle{definition}
\newtheorem{definition}[theorem] {Definition}

\theoremstyle{remark}
\newtheorem{remark}[theorem] {Remark}

\numberwithin{equation}{section}

\usepackage{amsfonts}
\usepackage{amssymb}
\usepackage{amscd}

\newcommand{\R}{{\mathbb R}}
\newcommand{\Z}{{\mathbb Z}}
\newcommand{\N}{{\mathbb N}}

\newcommand{\PP}{{\mathcal P}}
\newcommand{\F}{{\mathcal F}}

\newcommand{\CC}{{\mathbb C}}

\newcommand{\E}{{\mathcal E}}
\newcommand{\al}{{\alpha}}
\newcommand{\la}{{\lambda}}
\newcommand{\sa}{{\sigma}}

\newcommand{\iy}{{\infty}}

\newcommand{\vep}{{\varepsilon}}
\newcommand{\g}{{\gamma}}
\newcommand{\de}{{\delta}}

\newcommand{\be}{{\beta}}

\newcommand{\bna}{\begin{eqnarray}}
\newcommand{\ena}{\end{eqnarray}}
\newcommand{\ba}{\begin{eqnarray*}}
\newcommand{\ea}{\end{eqnarray*}}
\newcommand{\beq}{\begin{equation}}
\newcommand{\eeq}{\end{equation}}
\DeclareMathOperator*{\esssup}{ess\,sup}

\begin{document}

\title[Constants in  Markov--Bernstein--Nikolskii Inequalities]
{Asymptotics of Sharp Constants in Markov--Bernstein--Nikolskii
 Type Inequalities with Exponential Weights}
\author{Michael I. Ganzburg}
 \address{212 Woodburn Drive\\ Hampton,
 VA 23668\\USA}
 \email{michael.ganzburg@gmail.com}
 \keywords{Exponential weights, Markov--Bernstein inequality, algebraic polynomials,
 entire functions of exponential type}
 \subjclass[2010]{41A17, 41A10}
 \begin{abstract}
  We prove that the sharp constant
  in the univariate Bernstein--Nikolskii inequality
  for entire functions of exponential type is the limit of the
  sharp constant in the V. A. Markov type inequality with an exponential weight
   for coefficients of an algebraic polynomials of degree n as $n\to\iy$.
 \end{abstract}
 \maketitle

 \section{Introduction}\label{S1}
 \noindent
\setcounter{equation}{0}
In this paper we discuss limit relations between the sharp constants in the univariate
V. A. Markov--Bernstein--Nikolskii type
inequalities  with exponential weights
for algebraic polynomials and entire functions of exponential type.\vspace{.1in}\\
\textbf{Notation and Preliminaries.}
Throughout the paper $C,\,C_1,\,C_2,\ldots$ denote positive
constants independent
of essential parameters.
 Occasionally we indicate dependence on or independence of certain parameters.
 The same symbol $C$ does not
 necessarily denote the same constant in different occurrences.

Let
$\N:=\{1,\,2,\,\ldots\},\,\Z_+:=\{0,\,1,\,\ldots\}$,
$\R$ be the set of all real numbers,
$\CC$ be the set of all complex numbers,
$\PP_n$ be the set
 of all algebraic polynomials with complex coefficients
 of degree at most $n,\,n\in \Z_+$,
 and $B_\sa$ be
 the set of all complex-valued entire
 functions of exponential type $\sa > 0$.

 Let $W:\Omega\to [0,\iy]$ be an integrable weight on a measurable
 subset $\Omega$ of $\R$, and let
 $L_{r,W}(\Omega)$ be a weighted space of all measurable complex-valued
  functions $F:\Omega\to \CC$
  with the finite quasinorm
 \ba
 \|F\|_{L_{r,W}(\Omega)}:=\left\{\begin{array}{ll}
 \left(\int_\Omega\vert F(x)W(x)\vert^r dx\right)^{1/r}, & 0<r<\iy,\\
 \esssup_{x\in \Omega} \vert F(x)\vert W(x), &r=\iy.
 \end{array}\right.
 \ea
 In the nonweighted case ($W =1$), we set
$$
\|\cdot \|_{L_r (\Omega)} : =\|\cdot \|_{L_{r,1} (\Omega)} ,
\quad \ L_r (\Omega) := L_{r,1}
(\Omega),\qquad 0<r\leq \infty .
$$
The quasinorm $\|\cdot\|_{L_{r,W}(\Omega)}$
allows the following "triangle" inequality:
\beq\label{E1.1}
\left\|\sum_{j=1}^l F_j\right\|^{\tilde{r}}_{L_{r,W}(\Omega)}
 \le \sum_{j=1}^l \left\|F_j\right\|^{\tilde{r}}_{L_{r,W}(\Omega)},
 \qquad F_j\in L_{r,W}(\Omega),\qquad
 1\le j\le l,
 \eeq
 where $l\in\N$ or $l=\iy$ and $\tilde{r}:=\min\{1,r\}$ for $r\in(0,\iy]$.

 Throughout the paper we assume that $W:I\to (0,\iy)$ is
 an exponential weight of the form
 $W(x)=\exp[-Q(x)]$ defined on a bounded or unbounded interval
  $I=(-c,c),\,0<c\le\iy$, where
  $Q$ is a continuous function on $I$.

  A function $F\colon (0,c)\to (0,\infty )$ is said to be quasi-increasing
if there exists a constant $C>0$ such that
$
F(x) \leq C F(y),\, 0<x\leq y <c.
$
The following definition describes the class of weights that
 we use in this paper
 (see \cite[Definition 1.1]{LL2001} and \cite[Definition 1.4.5]{G2008}).

 \begin{definition}\label{D1.1}
Let $W =e^{-Q}$, where $Q\colon I\to [0,\infty )$ satisfies
the following properties:

\begin{enumerate}
\item[{\rm (a)}] $Q$ is even in $I$ and $Q(0) =0$.
\item[{\rm (b)}] $Q^\prime$ is continuous in $I$.
\item[{\rm (c)}] $Q^{\prime\prime}(x)$ exists and $Q''(x)>0,\,x\in(0,c)$.
\item[{\rm (d)}] $\lim\limits_{x\to c-} Q(x) =\infty$.
\item[{\rm (e)}] The function
\beq\label{E1.2}
T(x): =\frac{x Q'(x)}{Q(x)}\, ,\qquad x\neq 0,
\eeq
is quasi-increasing in $(0,c)$ with
\begin{equation}\label{E1.3}
\Lambda :=\inf_{x\in (0,c)} T(x)>1.
\end{equation}
\item[{\rm (f)}] There exists a constant $C>0$ such that for all $x\in (0,c)$,
\ba
\frac{Q'' (x)}{|Q' (x)|} \leq C\frac{|Q' (x)|}{Q(x)}\, .
\ea
\end{enumerate}
Then we write $W\in \F (C^2)$.
\end{definition}
Note that properties (a), (b), and (c) of Definition \ref{D1.1} imply that
$Q$ is positive and increasing on $(0,c)$.
More classes of weights are discussed in
\cite[Sect. 1.2]{LL2001} and \cite[Sect.  1.4]{G2008}.

The behaviour of the function $T$ defined by \eqref{E1.2} divides
$\F(C^2)$ into two classes. In the case $I =\mathbb{R}$,
a weight $W\in \F(C^2)$,  satisfying
the condition
$
\sup_{x\in \mathbb{R}} T(x)<\infty ,
$
is called \emph{a Freud weight}. A typical example of such a weight is
\ba
W_\alpha (x): =\exp ( -|x|^\alpha ),\qquad \alpha >1, \quad I =\mathbb{R}.
\ea
A weight $W\in \F(C^2)$, satisfying the condition
$
\sup_{x\in\mathbb{R}} T(x) =\lim_{x\to\infty} T(x) =\infty ,
$
is called \emph{an Erd\"os weight}.
In particular, any weight $W\in\F(C^2)$ on a bounded interval
 $(-c,c)$ is an Erd\"os
weight.
 A typical example of an Erd\"os weight for the unbounded interval is
\ba
W_{\al,\ell}(x):=\exp \bigl( -\exp_\ell (|x|^\alpha ) +\exp_\ell (0)\bigr) ,
\qquad \alpha >1, \ \ \ell \geq 1, \ \ I =
\mathbb{R},
\ea
where $\exp_0 (x): =x$ and
$
\exp_k (x): = \exp \bigl( \exp_{k -1} (x)\bigr),\,1\leq k\leq \ell .
$

For $I=\R$ and a weight $W =e^{-Q}$,
$Q$ has at most polynomial growth on $\mathbb{R}$
if $W$ is a Freud weight, and
$Q$ has faster than polynomial growth on $\mathbb{R}$
if $W$ is an Erd\"os  weight.
These and other properties of Freud and Erd\"os weights along with
more examples can be found in \cite{LL2001} (see also \cite{G2008}).

Next, we define two constants that play an important role
 in orthogonal polynomials
for and weighted approximation with exponential weights.
Let $W=e^{-Q}\in\F(C^2)$ and
let $a_n=a_n(Q)\in(0,c)$ be the $n$th Mhaskar--Rakhmanov--Saff number defined
 as the positive
root of the equation
$$
n =\frac{2}{\pi} \int^1_0 \frac{a_n x Q'(a_n x)}{\sqrt{1-x^2}} \, dx,
\qquad n\in\N,
$$
(see \cite{MS1984,R1984,MS1985,LL2001}). Further, the number
\ba
b_n =b_n (Q) := \frac{2}{\pi} \int^1_0 \frac{Q' (a_nx)\sqrt{1 -x^2}}{x}\, dx
+\frac{n}{a_n},\qquad n\in\N,
\ea
was defined in \cite[Eq. 1.2.2]{G2008} in order to replace $n$ in
sharp constants of
nonweighted approximation theory.
Note that for  Freud weights,
\beq\label{E1.8}
n/a_n\leq b_n \leq (1 +C) n/ a_n,\qquad n\in\N,
\eeq
 and for Erd\"os weights,
 \beq\label{E1.9}
 b_n= (n/a_n) (1 + o(1))
\eeq
as $n\to\infty$.
We also note that
\beq\label{E1.9a}
\lim_{n\to\iy} b_n=\iy
\eeq
and
\beq\label{E1.9b}
b_n\le Cn,\qquad n\in\N.
\eeq
Relations \eqref{E1.8}, \eqref{E1.9}, and \eqref{E1.9a} were proved in
\cite[Proposition 3.2.2]{G2008}, while \eqref{E1.9b} follows from
\eqref{E1.8}, \eqref{E1.9}, and increasing of $a_n,\,n\in\N$
(see \cite[Lemma 2.13]{LL2001}).

For example (see \cite[Sects. 10.1 and 10.5]{G2008})
for the weight $W_\al$,
\ba
a_n=\left(\frac{2^{\al-2}\Gamma^2(\al/2)}{\Gamma(\al)}\right)^{1/\al}
n^{1/\al},
\qquad
b_n=\frac{\al n}{(\al-1) a_n}
=\frac{\al}{(\al-1)}\left(\frac{\Gamma(\al)}{2^{\al-2}\Gamma^2(\al/2)}
\right)^{1/\al} n^{1-1/\al},
\ea
where $\Gamma(z)$ is the gamma function,
and for the weight $W_{\al,\ell}$,
\ba
a_n=\left(\log_l(n)\right)^{1/\al}(1+o(1)),\qquad
b_n=\frac{n}{\left(\log_l(n)\right)^{1/\al}}(1+o(1)),\qquad
n\to\iy,
\ea
where
$\log_0(x):=x,\,\log_k(x)=\log\left(\log_{k-1}(x)\right),\,k\in\N.$
For a bounded interval, $a_n=c(1+o(1))$ and $b_n=(n/c)(1+o(1))$ as $n\to\iy$.
Similarly in the nonweighted case ($W=1$) and for the interval
$I=(-1,1)$, we assume that
$a_n=1+o(1)$ and $b_n=n(1+o(1))$ as $n\to\iy$.
\vspace{.12in}\\
\textbf{V. A. Markov--Bernstein--Nikolskii Type Inequalities.}
 Next, we define  sharp constants in univariate
 V. A. Markov--Bernstein--Nikolskii type
inequalities for algebraic polynomials
and entire functions of exponential type. Let
\bna
&&M_{p,N,n}(W):=b_n^{-N-1/p}
\sup_{P\in\PP_n\setminus\{0\}}\frac{\left\vert P^{(N)}(0)\right\vert}
{\|P\|_{L_{p,W}(I)}},\label{E1.10}\\
&& M_{p,N,n}^*(W):=b_n^{-N-1/p}
\sup_{P\in\PP_n\setminus\{0\}}\frac{\left\vert P^{(N)}(0)\right\vert}
{\|P\|_{L_{p,W}([-a_n,a_n])}},\label{E1.11}\\
&&E_{p,N}:=\sa^{-N-1/p}
\sup_{f\in (B_\sa\cap L_p(\R))\setminus\{0\}}
\frac{\left\|f^{(N)}\right\|_{L_\iy(\R)}}
{\|f\|_{L_p(\R)}}.\label{E1.12}
\ena
Here, $p\in(0,\iy],\,N\in\Z_+$, and $n\in\N$.
Note that $E_{p,N}$ is a nonweighted sharp constant, and it does not depend
on $\sa$ (see \cite{GT2017} for the proof). Therefore, we can assume that
$\sa=1$ in \eqref{E1.12}.
The exact values of $E_{p,N}$ are known only in the following cases
(see \cite[Sect. 1]{GT2017}):
\beq\label{E1.16}
E_{\iy,N}=1,\qquad E_{2,N}=(\pi(2N+1))^{-1/2},
\eeq
while the close estimates
$ 0.5409/\pi <E_{1,0}<0.5484/\pi$ were proved by Gorbachev \cite{G2005}.

The first sharp constant in the
nonweighted inequality for polynomial coefficients
 was found by V. A. Markov \cite{M1892} (see also
 \cite[Eq. (5.1.4.1)]{MMR1994})  in the form ($I=(-1,1)$ and $n\in\N$)
 \beq\label{E1.13}
 M_{\iy,N,n}(1)
 =n^{-N}
 \left\{\begin{array}{ll}
 \left\vert T_{n-1}^{(N)}(0)\right\vert, &n-N\,\mbox{is odd},\\
 \left\vert T_{n}^{(N)}(0)\right\vert, &n-N\,\mbox{is even}
 \end{array}\right.
 =(1+o(1))E_{\iy,N}
 \eeq
 as $n\to\iy$, where $T_n\in\PP_n$ is the Chebyshev polynomial of the first kind.
 Labelle \cite{L1969}  found $M_{2,N,n}(1)$
 for $I=(-1,1)$, and it turns out that
\beq\label{E1.14}
M_{2,N,n}(1)=(1+o(1))E_{2,N}
\eeq
as $n\to\iy$.
The author \cite[Theorem 1.1]{G2017}
 extended \eqref{E1.13} and \eqref{E1.14} to any $p\in(0,\iy]$
in the form
\beq\label{E1.15}
\lim_{n\to\iy}M_{p,N,n}(1)=E_{p,N},\qquad I=(-1,1).
\eeq
Multivariate versions of \eqref{E1.15} were obtained in
\cite{G2019, G2020a, G2020c}.

Certain properties of $M_{p,0,n}(W)$ for ultraspherical weights were
discussed by Arestov and Deikalova \cite{AD2015}. Some asymptotics for
sharp constants in V. A. Markov--Bernstein--Nikolskii type inequalities
with power and  ultraspherical weights were obtained by the author
\cite[Theorems 4.1 and 4.3]{G2019}.

In this paper we obtain a weighted version of relation \eqref{E1.15} for
exponential weights. The following estimates for $W\in \F(C^2)$
and more general weights were obtained by the author
\cite[Theorem 2.3.2 (b)]{G2008}
(see also Lemma \ref{L2.6}):
\beq\label{E1.16a}
M_{p,N,n}^*(W)\le M_{p,N,n}(W)\le C\sqrt{N+1}(1-\vep)^{-N},
\qquad \vep\in(0,1),\quad 0\le N\le n,\quad n\in\N.
\eeq
  \textbf{Main Results and Remarks.}
  Our major result discusses the limit relations between $M_{p,N,n}(W),
  \linebreak
  M_{p,N,n}^*(W)$, and $E_{p,N}$. In particular, we find
  an asymptotic behaviour of the sharp constants in inequality
  \eqref{E1.16a}.

  \begin{theorem} \label{T1.2}
 If $W\in\F(C^2),\,N\in\Z_+$, and $p\in(0,\iy]$, then
 \beq \label{E1.17}
  \lim_{n\to\iy}M_{p,N,n}(W)
  =\lim_{n\to\iy}M_{p,N,n}^*(W)=E_{p,N}.
 \eeq
 \end{theorem}
 Combining Theorem \ref{T1.2} with relations \eqref{E1.16}, we
 arrive at the following corollary:

 \begin{corollary}
 If $W\in\F(C^2)$ and $N\in\Z_+$, then
 \ba
  &&\lim_{n\to\iy}M_{\iy,N,n}(W)
  =\lim_{n\to\iy}M_{\iy,N,n}^*(W)=1,\\
  &&\lim_{n\to\iy}M_{2,N,n}(W)
  =\lim_{n\to\iy}M_{2,N,n}^*(W)=(\pi(2N+1))^{-1/2}.
 \ea
 \end{corollary}

 \begin{remark}\label{R1.4}
In definitions \eqref{E1.10},
\eqref{E1.11} and \eqref{E1.12} of the sharp constants we
discuss only complex-valued functions $P$ and $f$. We can define similarly
the "real" sharp constants if the suprema in \eqref{E1.10},
\eqref{E1.11} and \eqref{E1.12}
 are taken over all real-valued functions
on $\R$ from $\PP_{n}\setminus\{0\}$
and $(B_\sa\cap L_p(\R))\setminus\{0\}$, respectively.
It turns out that the "complex" and "real" sharp constants coincide.
 For $W=1,\,I=(-1,1)$, this fact was proved in \cite[Sect. 1]{G2017} (cf.
\cite[Theorem 1.1]{GT2017} and \cite[Remark 1.5]{G2020a}),
and the case of exponential weights can be proved similarly.
In addition, Theorem \ref{T1.2} is
 formulated for the "complex"  sharp constants.
 However, this result remains valid for the "real" ones as well.
 The proof of the real version of Theorem \ref{T1.2}
 does not change compared with the complex one
 if we take into account  Remark \ref{R2.5} from Section \ref{S2}.
\end{remark}

 \begin{remark}\label{R1.5}
 Theorem \ref{T1.2} shows that the sharp constants in the weighted
 $L_p$-inequalities for the $N$th coefficient of a polynomial are
 asymptotically equal to $E_{p,N}b_n^{N+1/p}/N!$, where by
 \eqref{E1.8} and \eqref{E1.9},
 $b_n\sim n/a_n$ and for Erd\"os weights
 $b_n=(n/a_n)(1+o(1))$,
 as $n\to\iy$.
 Note that the sharp dependence on $n$ of the sharp constant
 in the A. A. Markov--Nikolskii type inequality
 with an exponential weight $W\in\F(C^2)$
is supposed to be $\left((n/a_n)\sqrt{T(a_n)}\right)^{N+1/p}$
(see \cite[Corollary 10.2 and Theorem 10.3]{LL2001}) compared with
$\left(n/a_n\right)^{N+1/p}$ in Theorem \ref{T1.2}. An asymptotic for
this constant is unknown. However,
the asymptotic behaviour of
  the sharp constant in the
classical nonweighted inequality of different metrics
was found in \cite[Corollary 4.6]{G2019} in the following form:
\beq\label{E1.17a}
\lim_{n\to\iy}n^{-2/p}\sup_{P\in\PP_n\setminus\{0\}}\frac{\|P\|_{L_{\iy}([-1,1])}}
{\|P\|_{L_{p}([-1,1])}}
=2^{1/p}\sup_{f\in (B_1\cap L_{p,W^*}(\R))\setminus\{0\}}
\frac{\left\vert f(0)\right\vert}
{\|f\|_{L_{p,W^*}(\R)}},\qquad p\in [1,\iy),
\eeq
where $W^*(x):=\vert x\vert^{1/p}$. A different version of  \eqref{E1.17a} for
$p\in (0,\iy)$ was proved in \cite[Theorem 1.4]{G2017} (see also
\cite[p. 94]{G2017}).
 \end{remark}

 The proof of Theorem \ref {T1.2} is presented in
Section \ref{S3}.
It follows general ideas developed in \cite[Corollary 7.1]{G2020b}.
Section \ref{S2} contains certain properties of functions
from $B_\sa$ and polynomials from $\PP_{n}$.

 \section{Properties of Entire Functions  and Polynomials}\label{S2}
 \noindent
\setcounter{equation}{0}
To prove Theorem \ref{T1.2}, we need several lemmas about
 certain properties of entire functions of exponential type
   and algebraic
 polynomials.
 We start with known properties of entire functions
 of exponential type.
 \begin{lemma}\label{L2.1}
 (a) The following crude Bernstein and Nikolskii type inequalities hold true:
 \bna
 &&\left\|f^{(s)}\right\|_{L_{\iy}(\R)}
  \le C
  \left\|f\right\|_{L_{\iy}(\R)},\qquad f\in B_\sa\cap L_\iy(\R),
  \quad s\in Z_+,\label{E2.1a}\\
   &&\left\|f\right\|_{L_{\iy}(\R)}
  \le C
  \left\|f\right\|_{L_{p}(\R)},\qquad f\in B_\sa\cap L_p(\R),\quad
 p\in(0,\iy),\label{E2.1b}
  \ena
  where $C$ is independent of $f$.\\
  (b) If $f\in B_\sa\cap L_p(\R),\,p\in(0,\iy)$, then
  \beq\label{E2.1c}
  \lim_{\vert x\vert\to\iy} f(x)=0.
  \eeq
  \end{lemma}
  \proof
  The proofs of \eqref{E2.1a} and \eqref{E2.1b} can be found
  in \cite[Theorem 11.3.3]{B1954} and \cite[Eq. 4.9(29)]{T1963}, respectively.
  The proof of a multivariate version of statement (b),  given in
  \cite[Theorem 3.2.5]{N1969} for $p\in[1,\iy)$, is long and difficult.
   For the reader's
    convenience, we present a shorter and more elementary proof of (b)
    for $p\in(0,\iy)$.

    If \eqref{E2.1c} is not valid, then there exist $\vep>0$ and
     a number sequence
    $\{x_n\}_{n=1}^\iy$ such that $\lim_{n\to\iy}\vert x_n\vert=\iy$ and
    $\inf_{n\in\N}\vert f(x_n)\vert\ge \vep$.
     Without loss of generality we can assume that $0<x_1<x_2<\ldots$.
      Setting $x_{n_1}:=x_1$ and $y_0:=0$
      and recalling that $f\in L_p(\R)$ and $f$ is continuous on $\R$,
       we can construct by induction
      a subsequence  $\{x_{n_k}\}_{k=1}^\iy$ and a sequence
      $\{y_k\}_{k=1}^\iy$ such that
      \ba
      \lim_{k\to\iy} x_{n_k}=\iy,\,y_k>x_{n_k}>y_{k-1},\,
      \left\vert f\left(x_{n_k}\right)\right\vert\ge \vep,\,
      \left\vert f\left(y_k\right)\right\vert= \vep/2,\,
      \inf_{x\in\left[x_{n_k},y_k\right]}\vert f(x)\vert\ge \vep/2,\,k\in\N.
      \ea
      Next, setting $\la_k:=y_k-x_{n_k},\,k\in\N$, we obtain
      \ba
      (\vep/2)^p\sum_{k=1}^\iy\la_k
      \le \sum_{k=1}^\iy\int_{x_{n_k}}^{y_k}\vert f(x)\vert^pdx
      \le \|f\|^p_{L_p(\R)}.
      \ea
      Therefore, $\lim_{k\to\iy}\la_k=0$, while
      by \eqref{E2.1a} and \eqref{E2.1b},
      \ba
      \vep/2\le  \left\vert f\left(x_{n_k}\right)-f(y_k)\right\vert
      \le \left\|f^\prime\right\|_{L_\iy(\R)}\la_k
      \le C\|f\|_{L_p(\R)}\la_k,\qquad k\in\N.
      \ea
       This contradiction proves statement (b). \hfill $\Box$

  Next, we need the following version of the compactness
      theorem for entire functions of exponential type.

     \begin{lemma}\label{L2.2}
     Let $\E$ be the set  of all entire functions $f(z)=\sum_{k=0}^\iy c_kz^k$,
     satisfying the following condition: for any $\de>0$ there exists a constant $C(\de)$,
     independent of $f$ and $k$, such that
     \beq\label{E2.1d}
     \vert c_k\vert \le \frac{C(\de)(1+\de)^k}{k!},\qquad k\in\Z_+.
     \eeq
     Then for any sequence $\{f_n\}_{n=1}^\iy\subseteq\E$ there exist
     a subsequence $\{f_{n_m}\}_{m=1}^\iy$ and a function $f_0\in B_1$
     such that for every $s\in\Z_+,\,\,\lim_{m\to\iy}f_{n_m}^{(s)}=f_0^{(s)}$
     uniformly on each compact subset of $\CC$.
     \end{lemma}
The lemma was proved in \cite[Lemma 2.6]{G2017}.

 Further, we discuss  estimates of the error of weighted polynomial approximation
  for functions from $B_\tau$.

 \begin{lemma}\label{L2.3}
 Let $W\in\mathcal{F} (C^2)$. Then there exist numbers $\delta_1=\delta_1 (W)
\in\left( 0,\frac{2(\Lambda -1)}{3(\Lambda +1)}\right),\, \delta_2 =\delta_2
(W)>0$, and a constant $C_1=C_1 (W)>0$ such that for every $\tau \in
\bigl(0,b_n
\left(1 -C_1n^{-\delta_1} \right)\bigr]$,
any $g\in B_\tau\cap L_\iy(\R)$, and each $k\in\N$
there exists a polynomial $P_k\in\PP_k$ such that
the following estimate holds:
\begin{equation}\label{E2.1}
\left\|g-P_k\right\|_{L_{r,W}(I)}
 \leq C_2 k^\gamma\exp \left( -C_3\, k^{\delta_2}
\right)\|g\|_{L_\iy(\R)} ,\qquad 0<r\leq \infty .
\end{equation}
Here, $\Lambda$ is defined by \eqref{E1.3}, and
$C_2,\,C_3$, and $\gamma$ are positive constants independent of $k$ and $g$.
\end{lemma}

 This lemma follows from a more general result proved in
 \cite[Theorem 2.2.1]{G2008}.
 Lemma \ref{L2.3} is a weighted version of estimates obtained
 by Bernstein \cite{B1946} (see also \cite[Sect. 5.4.4]{T1963} and
  \cite[Appendix, Sect. 85]{A1965}). More precise and more general
  nonweighted inequalities were proved by the author
  in \cite{G1982} and \cite{G1991}.

  \begin{lemma}\label{L2.4}
  Let $W\in\mathcal{F} (C^2)$.
  Then for any $\tau\in \bigl(0,b_n
\left(1 -C_1n^{-\delta_1} \right)\bigr],\,
  g\in B_\tau\cap L_\iy(\R)$ with $\|g\|_{L_\iy(\R)}\le C_4$,
  and $n\in\N$,
  there exists a polynomial
  $P_n\in\PP_n$ such that for every $s\in\Z_+$,
  $r\in(0,\iy]$, and $\eta\ge 0$,
  \beq\label{E2.2}
  \lim_{n\to\iy}n^\eta\left\|g^{(s)}-P_n^{(s)}\right\|_{L_{r,W}(I)}=0.
  \eeq
  Here, $C_1$ and $\de_1$ are the constants from Lemma \ref{L2.3} and
  the constants $C_4$ and $\eta$ are independent of $n$.
  \end{lemma}
  \proof
  First of all, for $P\in \PP_k,\,k\in\N$, and $r\in(0,\iy]$
  we need the following crude Markov-type inequality:
  \beq\label{E2.3}
  \left\|P^\prime\right\|_{L_{r,W}(I)}
  \le C_5(r,W)\, k\,\|P\|_{L_{r,W}(I)}.
  \eeq
  This inequality immediately follows from the estimates
  \ba
  \left\|P^\prime\right\|_{L_{r,W}(I)}
  \le C\, (k/a_k)\sqrt{T(a_k)}\|P\|_{L_{r,W}(I)}
  \ea
  (see \cite[Corollary 10.2]{LL2001}) and
  $T(a_k)\le C a_k^2,\,k\in\N$ (see inequality (3.38) in \cite[Lemma 3.7]{LL2001}).
  Here, $T$ is defined by \eqref{E1.2}, and
$C$ are constants independent of $k$ and $P$.

Next, let $\{P_k\}_{k=1}^\iy$ be the sequence of
  polynomials from Lemma \ref{L2.3}. Then using triangle inequality \eqref{E1.1}
  and inequalities \eqref{E2.3} and \eqref{E2.1}, we obtain
  \ba
  &&n^{\eta\tilde{r}}\left\|g^{(s)}-P_n^{(s)}\right\|_{L_{r,W}(I)}^{\tilde{r}}
  \le n^{\eta\tilde{r}}\sum_{k=n}^\iy
  \left\|(P_k-P_{k+1})^{(s)}\right\|_{L_{r,W}(I)}^{\tilde{r}}\\
  &&\le C_5^{s\tilde{r}}n^{\eta\tilde{r}}\sum_{k=n}^\iy (k+1)^{s\tilde{r}}
      \left\|P_k-P_{k+1}\right\|_{L_{r,W}(I)}^{\tilde{r}}\\
  &&\le  C_5^{s\tilde{r}}n^{\eta\tilde{r}}\sum_{k=n}^\iy (k+1)^{s\tilde{r}}
      \left(\left\|g-P_{k}\right\|_{L_{r,W}(I)}^{\tilde{r}}
      +\left\|g-P_{k+1}\right\|_{L_{r,W}(I)}^{\tilde{r}}\right)\\
  &&\le 2 C_2^{\tilde{r}}\, C_5^{s\tilde{r}}n^{\eta\tilde{r}}
  \sum_{k=n}^\iy (k+1)^{(s+\g)\tilde{r}}
  \exp \left( -C_3\, \tilde{r} \,k^{\delta_2}\right)\|g\|_{L_\iy(\R)}^{\tilde{r}}\\
  &&\le C_6 n^{\eta\tilde{r}}\int_{n}^\iy y^{(s+\g)\tilde{r}}
   \exp \left( -C_3\, \tilde{r} \,y^{\delta_2}\right)\,dy\,\|g\|_{L_\iy(\R)}^{\tilde{r}}\\
   &&\le  C_7C^{\tilde{r}}_4 n^{(s+\g+\eta)\tilde{r}}
   \exp \left( -C_3\, \tilde{r} \,n^{\delta_2}\right),
  \ea
where $C_2,\,C_3$, and $\de_2$ are constants from Lemma \ref{L2.3} and
$C_6$ and $C_7$ are constants independent of $n$.
Thus \eqref{E2.2} is established. \hfill $\Box$

\begin{remark}\label{R2.5}
  Note that if $g$ is a real-valued entire function in
  Lemmas \ref{L2.1} and \ref{L2.2}, then polynomials $P_n,\,n\in\N$,
  can be chosen real-valued as well.
  \end{remark}
 We also need a weighted estimate for coefficients of a polynomial.

\begin{lemma}\label{L2.6}
Let $W\in\F(C^2)$ and $p\in(0,\iy]$. Then for every polynomial
 $P(x)=\sum_{k=0}^nc_kx^k$
and any $\vep\in(0,1)$, the following inequality holds true:
\beq\label{E2.4}
\vert c_k\vert\le C_8(\vep,p,W)\,\frac{b_n^{k+1/p}}{(1-\vep)^k\, k!}
\|P\|_{L_{p,W}([-a_n,a_n])},\qquad 0\le k\le n.
\eeq
\end{lemma}
\proof
The inequality
\beq\label{E2.5}
\vert c_k\vert\le C_9(\vep,p,W)\,\frac{\sqrt{k+1}\,b_n^{k+1/p}}{(1-\vep/2)^k\, k!}
\|P\|_{L_{p,W}([-a_n,a_n])},\qquad \vep\in(0,1),\quad 0\le k\le n,
\eeq
was proved in \cite[Theorem 2.3.2]{G2008} for more general weights
(see also \eqref{E1.16a}).
 Then \eqref{E2.4}
follows from \eqref{E2.5} and an elementary inequality
\ba
\frac{\sqrt{k+1}}{(1-\vep/2)^k}
\le \frac{C_{10}(\vep)}{(1-\vep)^k},\qquad k\ge 0.
\ea
\hfill $\Box$

\section{Proof of Theorems \ref{T1.2}}\label{S3}
 \noindent
\setcounter{equation}{0}
We first prove the inequality
\beq\label{E3.1}
E_{p,N}\le\liminf_{n\to\iy}M_{p,N,n}(W).
\eeq

Let $f$ be a function from $B_1\cap L_p(\R),\,p\in(0,\iy]$.
Then $f\in B_1\cap L_\iy(\R)$ by \eqref{E2.1b}, and
$f^{(N)}\in B_1\cap L_p(\R)$ by \eqref{E2.1a} and \eqref{E2.1b}.

Let us first discuss the case $p\in(0,\iy)$. Then by Lemma \ref{L2.1} (b),
there exists $x_0\in\R$
such that $\left\|f^{(N)}\right\|_{L_\iy(\R)}=\left\vert f^{(N)}(x_0)\right\vert$.
Without loss of generality we can assume that $x_0=0$.
Next, setting $\be_n:=b_n\left(1 -C_1n^{-\delta_1} \right)$, we see that
the function $g_n(x):=f(\be_nx)$ belongs to $B_{\be_n}\cap L_\iy(\R)$
with $\|g_n\|_{L_\iy(\R)}=\|f\|_{L_\iy(\R)},\,n\in\N$.
In addition, recall that
$W(0)=1$ (by Definition \ref{D1.1}),
and $b_n\le Cn,\,n\in\N$ (by \eqref{E1.9b}).
Therefore, by Lemma \ref{L2.4} for $r=\iy,\,s=N,\,\eta=0$
and $r=p,\,s=0,\,\eta=1/p$,
 there exists a sequence of polynomials
$\{P_n\}_{n=1}^\iy$ such that
\bna
&&\lim_{n\to\iy}\left\vert g^{(N)}_n(0)-P_n^{(N)}(0)\right\vert
\le \lim_{n\to\iy}
\left\|g^{(N)}_n-P_n^{(N)}\right\|_{L_{\iy,W}(I)}=0,\label{E3.2}\\
&&\lim_{n\to\iy}b_n^{1/p}\left\|g_n-P_n\right\|_{L_{p,W}(I)}=0.\label{E3.3}
\ena
Then using \eqref{E3.2} and \eqref{E1.9a}  and taking account of definition
 \eqref{E1.10}, we obtain
\bna\label{E3.4}
\left\|f^{(N)}\right\|_{L_\iy(\R)}
&=&\left\vert f^{(N)}(0)\right\vert
=\be_n^{-N}\left\vert g^{(N)}_n(0)\right\vert
=\liminf_{n\to\iy}b_n^{-N}\left\vert g^{(N)}_n(0)\right\vert\nonumber\\
 &\le&  \lim_{n\to\iy}b_n^{-N}
 \left\vert g^{(N)}_n(0)-P_n^{(N)}(0)\right\vert
 +\liminf_{n\to\iy}b_n^{-N}\left\vert P_n^{(N)}(0)\right\vert\nonumber\\
 &=&\liminf_{n\to\iy}b_n^{-N}\left\vert P_n^{(N)}(0)\right\vert
 \le \liminf_{n\to\iy}\left(M_{p,N,n}(W)
 b_n^{1/p}\left\| P_n\right\|_{L_{p,W}(I)}\right).
 \ena
Next, using triangle inequality \eqref{E1.1} and
\eqref{E3.3}, we have
\bna\label{E3.5}
 \limsup_{n\to\iy}b_n^{\tilde{p}/p}
 \left\| P_n\right\|_{L_{p,W}(I)}^{\tilde{p}}
 \le \limsup_{n\to\iy}b_n^{\tilde{p}/p}\left(\|g_n-P_n\|_{L_{p,W}(I)}^{\tilde{p}}
 +\|g_n\|_{L_{p,W}(I)}^{\tilde{p}}\right)
 \le\|f\|_{L_p(\R)}^{\tilde{p}}.
 \ena
Combining (3.4) with (3.5), we arrive at (3.1) for $p\in(0,\iy)$.

In the case $p=\iy$, for any $\vep>0$
 there exists $x_0\in\R$
such that $\left\|f^{(N)}\right\|_{L_\iy(\R)}
<(1+\vep)\left\vert f^{(N)}(x_0)\right\vert$.
Without loss of generality we can assume that $x_0=0$. Then similarly to
\eqref{E3.4} and \eqref{E3.5} we can obtain the inequality
\beq\label{E3.6}
\left\|f^{(N)}\right\|_{L_\iy(\R)}
<(1+\vep)\liminf_{n\to\iy}M_{\iy,N,n}(W)
 \|f\|_{L_\iy(\R)}.
 \eeq
 Finally letting $\vep\to 0+$ in \eqref{E3.6},
 we arrive at \eqref{E3.1} for $p=\iy$.
  This completes the proof of \eqref{E3.1}.

 Further, we will prove the
 inequality
 \beq\label{E3.7}
\limsup_{n\to\iy}M_{p,N,n}^*(W)\le E_{p,N}
\eeq
  by constructing a nontrivial function $f_0\in B_1\cap L_p(\R)$,
   such that
 \beq \label{E3.8}
  \limsup_{n\to\iy}M_{p,N,n}^*(W)\le\left\|f_0^{(N)}\right\|_{L_\iy(\R)}/
\left\|f_0\right\|_{L_p(\R)}
  \le E_{p,N}.
 \eeq
 Since
 \beq \label{E3.9}
  M_{p,N,n}(W)\le M_{p,N,n}^*(W),
  \eeq
 inequalities \eqref{E3.1} and \eqref{E3.7} imply  \eqref{E1.17}.
 It remains to construct a nontrivial function $f_0$,
satisfying \eqref{E3.8}.
We first note that
\beq \label{E3.10}
\inf_{n\in\N}M_{p,N,n}^*(W)\ge C_{11}(p,N,W).
\eeq
This inequality follows immediately from  \eqref{E3.9} and \eqref{E3.1}.
Let $P_n\in\PP_n,\,n\in\N$, be a polynomial, satisfying the equality
\beq \label{E3.11}
M_{p,N,n}^*(W)=\frac{\left\vert P_n^{(N)}(0)\right\vert}
{b_n^{N+1/p}\left\|P_n\right\|_{L_{p,W}([-a_n,a_n])}}.
\eeq
The existence of an extremal polynomial $P_n$ in \eqref{E3.11}
can be proved by the standard compactness argument
(cf. \cite[Proof of Theorem 1.5]{GT2017}).
Next, setting $Q_n(x):=P_n(x/b_n)=\sum_{k=1}^nc_kx^k$,
 we have from \eqref{E3.11} that
\beq \label{E3.11a}
M_{p,N,n}^*(W)=\frac{\left\vert Q_n^{(N)}(0)\right\vert}
{\|Q_n\|_{L_{p,W(\cdot/b_n)}([-a_nb_n,a_nb_n])}}.
\eeq
Without loss of generality we can assume that
\beq \label{E3.12}
\left\vert Q_n^{(N)}(0)\right\vert=1.
\eeq
Then it follows from \eqref{E3.11a}, \eqref{E3.12}, and \eqref{E3.10}
that
\beq \label{E3.13}
\left\|Q_n\right\|_{L_{p,W(\cdot/b_n)}([-a_nb_n,a_nb_n])}
={1}/{M_{p,N,n}^*(W)}
\le {1}/{C_{11}(p,N,W)}.
\eeq
Further,  it follows from
inequality \eqref{E2.4} of
 Lemma \ref{L2.6} for $P=P_n$ and estimate \eqref{E3.13}   that
for every $\vep\in(0,1)$ and any $k,\,0\le k\le n,\,n\in\N$,
the following relations hold true:
\beq \label{E3.14}
\left\vert c_k\right\vert
=\frac{\left\vert P_n^{(k)}(0)\right\vert}
{b_n^{k}k!}
\le \frac{C_8(\vep,p,W) b_n^{1/p}
\left\|P_n\right\|_{L_{p,W}([-a_n,a_n])}}{(1-\vep)^kk!}
\le \frac{C_{8}(\vep,p,W)}{C_{11}(p,N,W)(1-\vep)^{k}k!}.
\eeq
Inequality \eqref{E3.14} holds true for all $k\in\Z_+$
if we set $c_k=0$ for $k>n,\,n\in\N$.
Therefore,
the polynomials
$Q_{n},\,n\in\N$, satisfy condition
\eqref{E2.1d}
of Lemma \ref{L2.2} with $\de:=\vep/(1-\vep)$
and $C(\de)=C_8/C_{11}$.
Thus $Q_n$ belongs to a set $\E$ of Lemma \ref{L2.2}, $n\in\N$.
Let $\{n_l\}_{l=1}^\iy$ be a subsequence of $\N$ such that
\beq \label{E3.15}
\limsup_{n\to\iy}M_{p,N,n}^*(W)=\lim_{l\to\iy}M_{p,N,n_l}^*(W).
\eeq
Then the polynomial sequence
$\{Q_{n_l}\}_{l=1}^\iy\subseteq \E$ satisfies all the conditions
of Lemma \ref{L2.2}.
Therefore, there exist
 a function $f_0\in B_1$ and a subsequence
 $\{Q_{n_{l_m}}\}_{m=1}^\iy$ such that
 \beq \label{E3.16}
 \lim_{m\to\iy}Q_{n_{l_m}}^{(s)}(x)=f_0^{(s)}(x),\qquad 0\le s\le N,
 \eeq
 uniformly on any interval $[-A,A],\,A>0$. Moreover, by
  \eqref{E3.12} and \eqref{E3.16},
  \beq \label{E3.17}
  \left\vert f_{0}^{(N)}(0)\right\vert=1.
\eeq
We also need the following relations:
\beq \label{E3.18}
\lim_{n\to\iy} a_nb_n=\iy,\qquad
\lim_{n\to\iy}\max_{x\in [-A,A]}\left(1-W(x/b_n)\right)=0.
\eeq
The first relation in \eqref{E3.18} follows from \eqref{E1.8}
and \eqref{E1.9}, while the second one follows from \eqref{E1.9a}
for every fixed $A>0$.

Next,
using consecutively  triangle inequality \eqref{E1.1} and relations
\eqref{E3.16}, \eqref{E3.18}, \eqref{E3.11a},  \eqref{E3.12},
and \eqref{E3.10},
we obtain for any interval $[-A,A],\,A>0$,
\bna \label{E3.19}
\|f_0\|_{L_p\left([-A,A]\right)}^{\tilde{p}}
&\le& \limsup_{m\to\iy}
\left(\|f_0-Q_{n_{l_m}}\|_{L_p\left([-A,A]\right)}^{\tilde{p}}
+\|Q_{n_{l_m}}\|_{L_p\left([-A,A]\right)}^{\tilde{p}}\right)\nonumber\\
&=& \limsup_{m\to\iy}\|Q_{n_{l_m}}\|_
{L_{p,W\left(\cdot/b_{n_{l_m}}\right)}\left([-A,A]\right)}^{\tilde{p}}\nonumber\\
&\le& \lim_{m\to\iy}\|Q_{n_{l_m}}\|_
{L_{p,W\left(\cdot/b_{n_{l_m}}\right)}
\left(\left[-a_{n_{l_m}}b_{n_{l_m}},
a_{n_{l_m}}b_{n_{l_m}}\right]\right)}^{\tilde{p}}\nonumber\\
&=&1/
\lim_{m\to\iy}M_{p,N,n_{l_m}}^*(W)
\le 1/C_{11}.
\ena
Therefore, letting $A\to \iy$ in \eqref{E3.19},
we see that $f_0$ is a nontrivial function from $B_1\cap L_p(\R)$,
by \eqref{E3.19} and \eqref{E3.17}. Thus for any interval
$[-A,A],\,A>0$, we obtain consecutively from \eqref{E3.15}, \eqref{E3.12},
 \eqref{E3.11a}, \eqref{E3.18}, \eqref{E3.16}, and \eqref{E3.17}
\bna \label{E3.20}
\limsup_{n\to\iy}M_{p,N,n}^*(W)
&=& \lim_{m\to\iy}\|Q_{n_{l_m}}\|_
{L_{p,W\left(\cdot/b_{n_{l_m}}\right)}
\left(\left[-a_{n_{l_m}}b_{n_{l_m}},
a_{n_{l_m}}b_{n_{l_m}}\right]\right)}^{-1}
\nonumber\\
&\le&\lim_{m\to\iy}\|Q_{n_{l_m}}\|_
{L_{p,W\left(\cdot/b_{n_{l_m}}\right)}([-A,A])}^{-1}\nonumber\\
&=&\lim_{m\to\iy}\|Q_{n_{l_m}}\|_
{L_{p}([-A,A])}^{-1}\nonumber\\
&=&\left\vert f_0^{(N)}(0)\right\vert/
\|f_0\|_{L_p([-A,A])}
\le \left\| f_0^{(N)}\right\|_{L_\iy(\R)}/
\|f_0\|_{L_p([-A,A])}.
\ena
Finally, letting $A\to \iy$ in \eqref{E3.20}, we arrive at \eqref{E3.8}.
\hfill$\Box$


\begin{thebibliography}{99}

\bibitem{A1965} N. I. Akhiezer, Lectures on the Theory of
Approximation. 2nd ed., Nauka, Moscow, 1965 (in Russian).

\bibitem{AD2015} V. Arestov, M. Deikalova, Nikol'skii inequality between
 the uniform norm and $L_q$-norm with ultraspherical weight of
 algebraic polynomials on an interval, Comput. Methods Funct. Theory,
  \textbf{15}  (2015), 689--708.

\bibitem{B1946} S. N. Bernstein, On the best approximation of continuous
 functions on the whole real axis by entire functions of given degree, V,
 Dokl. Akad. Nauk SSSR \textbf{54} (1946), 479--482 (in Russian).

 \bibitem{B1954} R. P. Boas, Entire Functions, Academic Press, New York, 1954.

\bibitem{G1982} M. I. Ganzburg, Multidimensional limit theorems of the theory
of best polynomial approximations, Sibirsk. Mat. Zh., \textbf{23} (1982),
no. 3, 30--47 (in Russian); English transl. in
 Siberian Math. J., \textbf{23}(1983), no. 3, 316--331.

 \bibitem{G1991} M. I. Ganzburg,	Limit theorems for the best
 polynomial approximation in the $L_\iy$-metric, Ukrain. Mat. Zh.
 \textbf{23} (1991), no. 3,
  336--342 (in Russian); English transl. in Ukrainian Math. J. \textbf{23}
  (1991), 299--305.

  \bibitem{G2008}  M. I. Ganzburg, Limit Theorems of Polynomial
Approximation with Exponential Weights,  Mem. Amer. Math. Soc.
\textbf{192} (2008), no. 897, 161 pp.

\bibitem{G2017} M. I. Ganzburg, Sharp constants in V. A. Markov-Bernstein
  type inequalities of different metrics, J. Approx. Theory
   \textbf{215} (2017), 92--105.

   \bibitem{G2019} M. I. Ganzburg, Sharp constants of approximation theory.
  II. Invariance theorems and  certain
  multivariate inequalities of different metrics,
   Constr. Approx. \textbf{50} (2019), 543--577.

   \bibitem{G2020a} M. I. Ganzburg, Sharp constants of approximation theory.
  III. Certain polynomial  inequalities of different metrics on convex sets,
   J. Approx. Theory \textbf{252} (2020), doi:10.1016/j.jat.2019.105351.

    \bibitem{G2020b} M. I. Ganzburg, Sharp constants of approximation theory.
  IV. Asymptotic relations in general settings,
  submitted; arXiv:2002.10512.

  \bibitem{G2020c} M. I. Ganzburg, Sharp constants of approximation theory.
  V. An asymptotic equality
related to polynomials with given Newton polyhedra,
  submitted; arXiv:2007.08439.

  \bibitem{GT2017} M. I. Ganzburg, S. Yu. Tikhonov, On sharp constants in
  Bernstein-Nikolskii inequalities, Constr. Approx.
  \textbf{45} (2017), 449--466.

  \bibitem{G2005} D. V. Gorbachev, An integral problem of Konyagin and the
 $(C,L)$-constants of
 Nikolskii, Teor. Funk., Trudy IMM UrO RAN, \textbf{11}, no. 2 (2005),
 72--91 (in Russian); English transl.
 in Proc. Steklov Inst. Math., Function Theory, Suppl. \textbf{2} (2005),
  S117--S138.

  \bibitem{L1969} G. Labelle, Concerning polynomials on the unit interval,
  Proc. Amer. Math. Soc. \textbf{20} (1969), 321--326.

  \bibitem{LL2001} Levin, A. L. and Lubinsky, D. S.,  Orthogonal
Polynomials for Exponential Weights, Springer, New York, 2001.

 \bibitem{M1892} V. A. Markov, On Functions Deviating Least from Zero
  in a Given Interval, Izdat. Imp. Akad. Nauk, St. Petersburg, 1892
  (in Russian); German transl. in Math. Ann. \textbf{77} (1916), 213--258.

  \bibitem{MS1984} H. N. Mhaskar,  E. B. Saff, Extremal problems
for polynomials with exponential weights, {\it Trans. Amer. Math. Soc.},
{\bf 285} (1984), 203--234.

\bibitem{MS1985} H. N. Mhaskar,  E. B. Saff, Where does the sup-norm
of a weighted polynomial live? {\it Constr. Approx.}, {\bf 1} (1985), 71--91.

\bibitem{MMR1994} G. V. Milovanovi\'{c}, D. S. Mitrinovi\'{c}, Th. M. Rassias,
 Topics in Polynomials: Extremal Problems, Inequalities, Zeros,
 World Scientific, Singapore, 1994.

 \bibitem{N1969} S. M. Nikolskii, Approximation of Functions of Several Variables
 and Imbedding Theorems, Nauka, Moscow, 1969 (in Russian); English edition:
 Die Grundlehren der Mathematischen Wissenschaften, Band 205, Springer-Verlag,
 New York-Heidelberg, 1975.

  \bibitem{R1984} E. A. Rakhmanov, On asymptotic properties of
polynomials orthogonal on the real axis, {\it Math. USSR Sbornik},
{\bf 47} (1984), 155--193.

 \bibitem{T1963} A. F. Timan, Theory of Approximation of Functions
  of a Real Variable, Fizmatgiz, Moscow, 1960 (in Russian);
  English edition: Pergamon Press, New York, 1963.






\end{thebibliography}
\end{document}